.
.
\font\script=eusm10.
\font\sets=msbm10.
\font\stampatello=cmcsc10.
\font\symbols=msam10.

\def\1{{\bf 1}}
\def\sgn{{\rm sgn}}

\def\square{\hbox{\vrule\vbox{\hrule\phantom{s}\hrule}\vrule}}
\def\defineq{\buildrel{def}\over{=}}
\def\defin{\buildrel{def}\over{\Longleftrightarrow}}
\def\doublesum{\mathop{\sum\sum}}
\def\integrale{\mathop{\int}}
\def\supporto{{\rm supp}\,}

\def\N{\hbox{\sets N}}
\def\R{\hbox{\sets R}}
\def\Z{\hbox{\sets Z}}
\def\Corr{\hbox{\script C}}
\def\EssBdd{\hbox{\symbols n}\,}

\par
\noindent
\centerline{\bf ON THE SYMMETRY OF ARITHMETICAL FUNCTIONS}
\centerline{\bf IN ALMOST ALL SHORT INTERVALS, V}
\medskip
\centerline{by}
\smallskip
\centerline{G.Coppola}
\bigskip
{
\font\eightrm=cmr8
\eightrm {
\par
{\bf Abstract.} We study the symmetry in short intervals of arithmetic functions with non-negative exponential sums.
}
\footnote{}{\par \noindent {\it Mathematics Subject Classification} $(2000) : 11N37, 11N25.$}
}
\bigskip
\par
\noindent {\bf 1. Introduction and statement of the results.}
\smallskip
\par
\noindent
We pursue the study of the symmetry in (almost all) short intervals of arithmetical functions $f$ (see [C1]), where this time we give (non-trivial) results for a new class of such (real) $f$; the key-property they have is a non-negative exponential sum (see Lemma 2), which is something we require, in order to get a kind of \lq \lq majorant principle\rq \rq; this allows us to \lq \lq smooth\rq \rq \thinspace our $f$ into a \lq \lq restricted\rq \rq \thinspace divisor function (see the Theorem), for which (see the Corollary) we apply non-trivial results (both from [C-S], on Acta Arithmetica, and [C2]).
\par
We need the following definitions. Here and in the sequel $h\to \infty$ and $h=o(N)$, if  $N\to \infty$.
$$
\Corr_{f_1,f_2}(a)\defineq \sum_{|a|<n\le N-|a|}f_1(n)f_2(n-a)\quad \forall a\in [-2h,2h]\cap \Z
$$
\par
\noindent
(are the \lq \lq {\stampatello mixed\rq \rq \thinspace correlations} of arithmetic, real, $f_1,f_2:\N \rightarrow \R$) {\stampatello say} $\sgn(r)\defineq {{|r|}\over r}$ $\forall r\in \R^{*}$, $\sgn(0)\defineq 0$: 
$$
I_{f_1,f_2}(N,h)\defineq \int_{h}^{N}\Big( \sum_{|n-x|\le h}f_1(n)\sgn(n-x)\Big) \Big( \sum_{|m-x|\le h}f_2(m)\sgn(m-x)\Big)dx
$$
\par
\noindent
(are the \lq \lq {\stampatello mixed\rq \rq \thinspace symmetry integral} of the arithmetical, real, functions $f_1$ \& $f_2$)
$$
I_f(N,h)\defineq \int_{h}^{N}\Big| \sum_{|n-x|\le h}f(n)\sgn(n-x)\Big|^2 dx
$$
\par
\noindent
(is the {\stampatello symmetry integral} of the arithmetical real function $f:\N \rightarrow \R$)
\smallskip
\par
\noindent $\underline{\hbox{\stampatello Remark}}$ Here $I_{f_1,f_2}$ \& $I_f$ (like, also, mixed correlations) depend only on $f,f_1,f_2$ values in $[1,N+h-1]\cap \N$; \indent \indent for all of these quantities it is essential to assume $h=o(N)$ (if $N\to \infty$), to avoid trivialities.
\par
\noindent
Write the {\stampatello divisor function} \enspace $d(n)\defineq \sum_{q|n}1$ \enspace {\stampatello and} \enspace $d_Q(n)\defineq \sum_{q|n,q\le Q}1$ \enspace the {\stampatello \lq \lq restricted\rq \rq \thinspace divisor function} (correspondingly $I_{f,d_Q}$ and $I_{d_Q}$). \enspace Also, $\ast$ is Dirichlet product (see the following) and $\mu$ is M\"obius function [T]; we'll indicate \lq \lq $\supporto$\rq \rq \thinspace for the support of our functions and abbreviate 
$$
A\EssBdd B \enspace \defin \enspace \forall \varepsilon>0 \enspace A\ll_{\varepsilon} N^{\varepsilon}B
$$
\par
\noindent
whence, for example, $d(n)\EssBdd 1$ (i.e., {\stampatello say}, the divisor function is {\stampatello \lq \lq essentially bounded\rq \rq}).
\par
\noindent
We define (see Lemma 1) \quad $W(a)\defineq \max(2h-3|a|,|a|-2h)$ \enspace $\forall a\in \Z \cap [-2h,2h]$, \enspace $\supporto(W)\subset [-2h,2h]$.
\medskip
\par
Then, we come to our main result, which we refer to as a {\stampatello \lq \lq majorant principle\rq \rq \thinspace } for $I_f$:
\smallskip
\par
\noindent $\underline{\hbox{\stampatello Theorem}}$. {\it Let } $h\to \infty$ {\it and } $h=o(N)$ {\it when } $N\to \infty$. {\it Assume that the arithmetical function } $f:\N \rightarrow \R$ ({\it which might, and actually will, depend on } $N$ {\it and } $h$) {\it has } $\supporto(f\ast \mu)\subset [1,Q]$, {\it together with the property}
$$
\forall \varepsilon>0 \quad f(n)\ll_{\varepsilon} N^{\varepsilon} \enspace \forall n\ll N \enspace (\hbox{\stampatello uniformly\thinspace in} \enspace 0<n\ll N)
$$
\par
\noindent
{\it and indicate with } $f(0)>0$ \enspace {\it a constant which may depend on } $N,h$ ({\stampatello i.e.} $f(0)=f_{N,h}(0)>0$). {\it Then}
$$
S_f(\alpha)\defineq \Re\sum_{0\le n\le N}f(n)e(n\alpha)\ge 0 \enspace \forall \alpha\in [0,1] 
\enspace \Rightarrow \enspace 
I_f(N,h)\EssBdd I_{f,d_Q}(N,h)+h^3+f(0)h^2+Qhf(0)+Qh^2.
$$
\medskip
\par
\noindent
An immediate consequence is:
\smallskip
\par
\noindent $\underline{\hbox{\stampatello Corollary}}$. {\it In the Theorem hypotheses}, {\stampatello assuming also } $h^2\ll Q$, \enspace $\theta \defineq {{\log h}\over {\log N}}<{1\over 2}$,
$$
S_f(\alpha)\ge 0 \enspace \forall \alpha\in [0,1] 
\Rightarrow 
I_f(N,h)\EssBdd Nh\sqrt{h}+Qhf(0)+Qh^2.
$$
\medskip
\par
\noindent
In fact, applying Cauchy-Schwarz inequality \enspace $I_{f,d_Q}(N,h)\le \sqrt{I_f(N,h)}\sqrt{I_{d_Q}(N,h)}$\enspace and
$$
f(n)\EssBdd 1 \enspace \Rightarrow \enspace I_f(N,h)\EssBdd Nh^2 \quad \hbox{\rm is\thinspace the\thinspace trivial\thinspace estimate, while}
$$
$$
\theta < {1\over 2} \enspace \Rightarrow \enspace I_{d_Q}(N,h)\EssBdd Nh, \quad \hbox{\stampatello see [C-S] \& compare [C2].}
$$
\par
\noindent
All of this gives \enspace $I_{f,d_Q}(N,h)\EssBdd Nh^{3/2}$, which is bigger than $\EssBdd h^3$ \enspace (due to \enspace $\theta<{1\over 2}<{2\over 3}$, here).
\bigskip
\par
\noindent $\underline{\hbox{\stampatello Remark}}$ The feature, both in the Theorem and in the Corollary,
$$
S_f(\alpha)\ge 0 \enspace \thinspace \enspace \forall \alpha \in [0,1]
$$
\par
\noindent
is too strong, as the requirement might be milder (in order to apply Lemma 2, see following section): 
$$
S_f\left({j\over q}\right)\ge 0, \enspace \forall j\in \Z_{q}^{*} \enspace \hbox{\stampatello and} \enspace \forall q\le Q
$$
\par
\noindent
(being \enspace $\Z_{q}^{*}$\enspace the reduced residue classes modulo $q$, here) suffices for our Theorem (and our Corollary).
\bigskip
\par
We have used and will use the notation, for {\stampatello Dirichlet product}
$$
\left( f_1 \ast f_2 \right)(n)\defineq \sum_{d|n}f_1(d)f_2(n/d) = \sum_{d|n}f_1(n/d)f_2(d) \quad \forall n\in \N;
$$
\par
\noindent
then, {\stampatello M\"obius inversion formula} reads \thinspace $f=g\ast \1 \Leftrightarrow g=f\ast \mu$, see [T] (whence $f\EssBdd 1 \Leftrightarrow g\EssBdd 1$).
\par
In the sequel \enspace $(j,q)=1$ \enspace indicates, as usual, that $j,q$ are coprime (no common prime divisors) and we'll write \enspace $j(\bmod \,\, q)$ \enspace for the residue classes modulo $q$ (doesn't matter if \thinspace $j\le q$ \thinspace or \thinspace $0\le |j|\le q/2$, here).
\par
\noindent
Also, we' ll follow the standard notation \enspace $e(n\alpha)\defineq e^{2\pi in\alpha}$ ($\forall n\in \N \; \forall \alpha \in \R$) \enspace for additive characters.

\bigskip

\par
\noindent $\underline{\hbox{\stampatello Acknowledgement}}$. The author thanks Professor Alberto Perelli for friendly and useful remarks during long and  enlightning conversations.

\bigskip

\par
We'll start with three Lemmas which (resp.ly) manage correlations, link mixed symmetry integrals with mixed correlations and, finally, give the essence of our, say, {\stampatello \lq \lq majorant principle\rq \rq.}

\medskip

\par
\noindent
The paper is organized as follows:
\smallskip
\item{---} in section 2 we state and prove our Lemmas;
\smallskip
\item{---} in section 3 we prove our Theorem;
\smallskip
\item{---} in section 4 we give an example.

\vfil
\eject

\par
\noindent {\bf 2. Lemmas.}
\smallskip
\par
\noindent $\underline{\hbox{\stampatello Lemma 0}}$. {\it Let } $h\to \infty$ {\it and } $h=o(N)$ {\it when } $N\to \infty$. {\it Assume } $a\in \Z$ {\it with } $0<|a|\le 2h$ {\it and}
$$
f_1:\N \rightarrow \R,\, f_2:\N \rightarrow \R \enspace \hbox{\stampatello are\thinspace such\thinspace that} \enspace \forall \varepsilon>0 \enspace f_1(n),f_2(n)\ll_{\varepsilon} N^{\varepsilon}\, \hbox{\stampatello uniformly}\, \forall n\ll N.
$$
\par
\noindent
{\it Then } $\forall \varepsilon>0$
$$
\qquad \thinspace \enspace \thinspace \qquad
\sum_{{2h<n<N-h}\atop {2h<n-a<N-h}}f_1(n)f_2(n-a)=\Corr_{f_1,f_2}(a)+{\cal O}_{\varepsilon}\left( N^{\varepsilon}h\right) \enspace \hbox{\stampatello and} \thinspace \enspace \thinspace \Corr_{f_1,f_2}(-a)=\Corr_{f_2,f_1}(a)+{\cal O}_{\varepsilon}\left( N^{\varepsilon}h\right).
$$
\par
\noindent $\underline{\hbox{\stampatello proof}}.\!$ Assume \enspace $a>0$\enspace : LHS becomes (recall \enspace $|a|\ll h$ \enspace and \enspace $f_1,f_2\ll_{\varepsilon} N^{\varepsilon}$, here)
$$
\sum_{2h+a<n<N-h}f_1(n)f_2(n-a)=\sum_{a<n\le N-a}f_1(n)f_2(n-a)+{\cal O}_{\varepsilon}\left( N^{\varepsilon}h\right)=\Corr_{f_1,f_2}(a)+{\cal O}_{\varepsilon}\left( N^{\varepsilon}h\right); 
$$
\par
instead, in the range \enspace $a<0$, $|a|\ll h$ (actually, $-2h\le a<0$), 
$$
\sum_{2h<n<N-h+a}f_1(n)f_2(n-a)=\sum_{-a<n\le N+a}f_1(n)f_2(n-a)+{\cal O}_{\varepsilon}\left( N^{\varepsilon}h\right)=\Corr_{f_1,f_2}(a)+{\cal O}_{\varepsilon}\left( N^{\varepsilon}h\right). 
$$
$$
\Corr_{f_1,f_2}(-a)=\sum_{|a|+a<m\le N-|a|+a}f_1(m-a)f_2(m)=\Corr_{f_2,f_1}(a)+{\cal O}_{\varepsilon}\left( N^{\varepsilon}h\right).\enspace \square
\leqno{\indent\hbox{\stampatello Finally},}
$$

\bigskip

\par
\noindent $\underline{\hbox{\stampatello Lemma 1}}$. {\it Let } $h\to \infty$ {\it and } $h=o(N)$ {\it when } $N\to \infty$. {\it Assume that }
$$
f_1:\N \rightarrow \R \enspace \hbox{\stampatello and} \enspace f_2:\N \rightarrow \R \enspace \hbox{\stampatello satisfy} \enspace \forall \varepsilon>0 \enspace f_1(n),f_2(n)\ll_{\varepsilon} N^{\varepsilon}, \hbox{\stampatello uniformly} \enspace \forall n\ll N.
$$
\par
\noindent
$$
I_{f_1,f_2}(N,h)=\sum_a W(a)\Corr_{f_1,f_2}(a)+{\cal O}_{\varepsilon}\left( N^{\varepsilon}h^3\right).
\leqno{{\it Then}\enspace \forall \varepsilon>0}
$$
\par
\noindent $\underline{\hbox{\stampatello proof}}.\!$ From the definition, exchanging sums and integral, LHS is, {\stampatello say}, 
$$
\enspace \thinspace \enspace \thinspace \enspace 
\doublesum_{{n,m\le N+h-1}\atop {0\le |n-m|\le 2h}}f_1(n)f_2(m)\integrale_{{h<x<N}\atop {|x-n|\le h,|x-m|\le h}}\sgn(x-n)\sgn(x-m)dx=\doublesum_{{n,m\le N+h-1}\atop {0\le |n-m|\le 2h}}f_1(n)f_2(m)\hbox{\script I}_{N,h}(m,n),
$$
\par
since \thinspace $x=h$ and $x=N$ \thinspace have no importance ($0-$measure) in the integral and \enspace $|x-n|\le h$, $|x-m|\le h$
\par
give \thinspace $|n-m|\le 2h$ (from triangle inequality); here the condition $h<x<N$ can be dispensed with into
$$
\doublesum_{{2h<n,m<N-h}\atop {0\le |n-m|\le 2h}}f_1(n)f_2(m)\hbox{\script I}_{N,h}(m,n)=\doublesum_{{2h<n,m<N-h}\atop {0\le |n-m|\le 2h}}f_1(n)f_2(m)W(|n-m|),
\eqno{\rm where}
$$ 
$$
\integrale_{{|x-n|\le h}\atop {|x-m|\le h}}\sgn(x-n)\sgn(x-m)dx=\int_{\max(n-h,m-h)}^{\min(n+h,m+h)}\sgn(x-n)\sgn(x-m)dx=W(|n-m|)
$$
\par
({\stampatello say, in accordance with } $W$ {\stampatello definition}, see the above) and, {\stampatello in fact, this is the main term}
$$
\qquad
\sum_a \doublesum_{{2h<n,m<N-h}\atop {n-m=a}}f_1(n)f_2(m)W(a)=\sum_a W(a)\!\!\sum_{{2h<n<N-h}\atop {2h<n-a<N-h}}f_1(n)f_2(n-a)=\sum_a W(a)\Corr_{f_1,f_2}(a)+{\cal O}_{\varepsilon}\left( N^{\varepsilon}h^3\right) 
$$
\par
from Lemma 0, with a good remainder; like also the terms completing LHS, above expanded : 
$$
\left( \doublesum_{{n\le 2h,\thinspace m\le N+h-1}\atop {0\le |n-m|\le 2h}}f_1(n)f_2(m)+\doublesum_{{m\le N+h-1,\thinspace N-h\le n\le N+h-1}\atop {0\le |n-m|\le 2h}}f_2(m)f_1(n)\right)\hbox{\script I}_{N,h}(m,n)\ll_{\varepsilon} N^{\varepsilon}h^3.\enspace \square
$$

\bigskip

\par
\noindent
We'll write, in the sequel, \enspace $\widehat{W}(\beta)\defineq \sum_{a}W(a)e(a\beta)$\enspace $\forall \beta \in \R$, \enspace the \lq \lq Discrete Fourier Transform\rq \rq \enspace of our \thinspace $W$.
\medskip
\par
\noindent $\underline{\hbox{\stampatello Lemma 2}}$. {\it Let } $g:\N \rightarrow \R$ {\it be such that } $\forall \varepsilon>0$ $g(q)\ll_{\varepsilon} N^{\varepsilon}$ $\forall q\le Q$ {\it and } $Q\ll N$ ($Q\to \infty$) {\it if } $N\to \infty$. {\it Then} 
$$
S_f(\alpha)\ge 0 \; \forall \alpha\in [0,1] \Rightarrow \sum_{q\le Q}{{g(q)}\over q}\sum_{j(\!\!\bmod q)}\widehat{W}\left( {j\over q}\right)S_f\left( -{j\over q}\right)\EssBdd \sum_{q\le Q}{1\over q}\sum_{j(\!\!\bmod q)}\widehat{W}\left( {j\over q}\right)S_f\left( -{j\over q}\right).
$$
\medskip
\par
\noindent
The {\stampatello proof} is immediate and simply applies the hypotheses (recall $\widehat{W}\ge 0$, [(1), {\stampatello Lemma} 4, C1]).

\vfil
\eject

\par
\noindent {\bf 3. Proof of the Theorem.}
\smallskip
\par
\noindent
We start giving the following :
\par
\noindent
{\stampatello in our hypotheses on $h,N$, writing} (in {\stampatello Lemma} 1) $f_1=f$, $f_2=g\ast \1$, {\stampatello with} $M:=\max \supporto(g)$,
$$
I_{f,g\ast \1}(N,h)=\sum_{a}W(a)\sum_{q\le M}g(q)\sum_{{0\le |n|\le N}\atop {n\equiv a(q)}}f(n)+{\cal O}_{\varepsilon}\left( N^{\varepsilon}\!\left( f(0)h^2 + Mhf(0)+Mh^2+h^3\right)\!\right)
\leqno{(*)}
$$
\par
\noindent
{\stampatello for} $g,f:\N \rightarrow \R$ {\stampatello any essentially bounded functions}. In fact, from Lemma 1 (see also Lemma 0)
$$
I_{f,g\ast \1}(N,h)=\sum_{a}W(a)\Corr_{f,g\ast \1}(a)+{\cal O}_{\varepsilon}\left( N^{\varepsilon}h^3\right)=\sum_{a}W(a)\sum_{q\le M}g(q)\sum_{{|a|<n\le N-|a|}\atop {n\equiv a(q)}}f(n)+{\cal O}_{\varepsilon}\left( N^{\varepsilon}h^3\right),
$$
\par
\noindent
which is, applying $h=o(N)$, here:
$$
\sum_{a}W(a)\sum_{q\le M}g(q)\sum_{{0\le n\le N}\atop {n\equiv a(q)}}f(n)+{\cal O}_{\varepsilon}\left( N^{\varepsilon}h^3\right)\!-\!\left(\sum_{q\le M}g(q)\!\!\right)\!\!\!\left(\sum_{a>0}W(a)f(a)\!\!\right)\!-\!f(0)\sum_{a}W(a)\sum_{q|a}g(q) ,
$$
\par
\noindent
whence $(\ast)$. {\stampatello Then, {\rm we} use} $(\ast)$ {\stampatello twice, one time for} $g=f\ast \mu$ (to treat $I_f$) {\stampatello and, soon after Lemma 2, with} $g=\1$ (supported in $[1,Q]$, to get $I_{f,d_Q}$).
 (In what  follows, we use \thinspace $\widehat{W}$, see before Lemma 2.)
\par
{\stampatello In fact, \enspace $g=f\ast \mu$ $\Rightarrow$ $I_f=I_{f,g\ast \1}$ \enspace and} \enspace $\supporto(g)\subset [1,Q]$ (see our Theorem hypotheses): 
$$
I_f(N,h)\EssBdd \left|\sum_{q\le Q}{{g(q)}\over q}\sum_{j(\!\!\bmod q)}\widehat{W}\left( {j\over q}\right)S_f\left( -{j\over q}\right)\right|+h^3+f(0)h^2+Qhf(0)+Qh^2
$$
$$
\EssBdd \sum_{q\le Q}{1\over q}\sum_{j(\!\!\bmod q)}\widehat{W}\left( {j\over q}\right)S_f\left( -{j\over q}\right)+h^3+f(0)h^2+Qhf(0)+Qh^2,
$$
\par
\noindent
from $(\ast)$ above, the orthogonality [D] of additive characters [V] and Lemma 2, with, {\stampatello say},
$$
S_f(\alpha):=\sum_{0\le n\le N}f(n)e(n\alpha),
$$
\par
\noindent
but we restrict to its real part (as defined in our theorem), {\stampatello since} 
$$
W\; \hbox{\stampatello even} \enspace \Rightarrow \enspace \widehat{W}\; \hbox{\stampatello even}.
$$
\par
\noindent
{\stampatello Finally, apply} $g=\1$ ({\stampatello with $\supporto(g)\subset [1,Q]$, here}) {\stampatello into} $(\ast).\enspace \square$
\bigskip
\par
\noindent {\bf 4. A first non-trivial (though non-optimal) application.}
\smallskip
\par
\noindent
From the remark above, we may restrict \thinspace $\alpha$ \thinspace to \enspace $\hbox{\script F}_Q := \{ j/q : j\le q, (j,q)=1, q\le Q\}$ ({\stampatello Farey fractions})
$$
S_f(\alpha)\ge 0 \enspace \forall \alpha \in \hbox{\script F}_Q 
\enspace \Rightarrow \enspace 
I_f(N,h)\EssBdd I_{f,d_Q}(N,h) + Qf(0)h + Qh^2
\leqno{(*)}
$$
\par
\noindent
and (see last remainder) we'll assume $\lambda \defineq {{\log Q}\over {\log N}}<1$ in our future papers applying the majorant principle.
\par
We want to get the condition $S_f\ge 0$ (see $(\ast)$ above) and we are turning upside-down the usual approach: instead of taking a class of functions $f$ to study, we start from the requirement $S_f\ge 0$ (at least on $\hbox{\script F}_Q$, see above) and want to understand for which class of functions $f:\N \rightarrow \R$ does it hold.
\par
The idea behind \lq \lq {\stampatello making $S_f$ positive}\rq \rq \thinspace is very easy: we think of the {\stampatello values} $f(n)$, $\forall n\neq 0$, as {\stampatello fixed}, {\stampatello choosing the mean-value} $f(0)$ in order {\stampatello to render} the whole sum $S_f$ {\stampatello positive} (actually, non-negative).
\par
The problem, now, is a definite one : {\stampatello choose $f(0)$ in order to ensure both $S_f\ge 0$ and that $Qf(0)h$ is a non-trivial remainder} (namely, \thinspace $Qf(0)h\ll Nh^2N^{-\delta}$, for some \thinspace $\delta>0$, fixed).
\par
Let's start from $S_f\ge 0$ on $\hbox{\script F}_Q$ : see that, $\forall n\neq 0$, \enspace $f(n)=\sum_{d|n}g(d)$, $\supporto(g)\subset [1,Q]$ gives ($f$ {\stampatello even})
$$
\sum_{0<|n|\le N}f(n)e(n\alpha) = \sum_{d\le Q}g(d)\sum_{0<|m|\le {N\over d}}e(d\alpha m):= \Sigma_1(\alpha)+\Sigma_2(\alpha)
$$
\par
\noindent
and we distinguish (as a kind of \lq \lq {\stampatello major arcs} \rq \rq \thinspace \& \lq \lq {\stampatello minor arcs}\rq \rq, here), say, with $\alpha={j\over q}$ :
$$
\Sigma_1(\alpha)\defineq \sum_{{d\le Q}\atop {d\equiv 0(\!\!\bmod q)}}g(d)\sum_{0<|m|\le {N\over d}}1 
\quad \& \quad
\Sigma_2(\alpha)\defineq \sum_{{d\le Q}\atop {d\not \equiv 0(\!\!\bmod q)}}g(d)\sum_{0<|m|\le {N\over d}}e(d\alpha m), 
$$
\par
\noindent
where we want an {\stampatello uniform} estimate over $\alpha \in \hbox{\script F}_Q$.
\par
Fix $\alpha \in \hbox{\script F}_Q$. \enspace Let $\Vert \alpha \Vert \defineq \min_{n\in \Z}|\alpha-n|$, $\forall \alpha \in \R$\enspace is the {\stampatello distance} of the real number \thinspace $\alpha$ \thinspace {\stampatello from the integers}. Then, \quad $g\ge 0 \Rightarrow \Sigma_1(\alpha)\ge 0$ \quad and (keeping the hypothesis \enspace $f\ast \mu \defineq g\ge 0$, now on) we don't care about $\Sigma_1(\alpha)$ value : it's a positive contribute to \enspace $S_f$, doesn't matter how much.
\par
\noindent
As regards the other sum, recall the hypothesis \thinspace $f\EssBdd 1$ \thinspace ($f$ \thinspace is {\stampatello essentially bounded}, on non-zero integers) entails (from $f,g$ definition, see above) \enspace $f\EssBdd 1 \enspace \Rightarrow \enspace g\EssBdd 1$, whence (see [D], chap.26) \enspace ${\displaystyle \sum_{0<|m|\le {N\over d}}e(d\alpha m)\ll {1\over {\left \Vert \alpha d\right \Vert}}}$
\par
\noindent
($\alpha \in \hbox{\script F}_Q$ $\Rightarrow$ $\alpha d$ is not an integer, $\forall d$ not a multiple of $q$) gives 
$$
\Sigma_2(\alpha)\ll \sum_{{d\le Q}\atop {d\not \equiv 0(\!\!\bmod q)}}g(d){1\over {\left \Vert {{jd}\over q}\right \Vert}}
\EssBdd \sum_{0<|r|\le {q\over 2}}{q\over {|r|}}\sum_{{d\le Q}\atop {d\equiv r\overline{j}(\!\!\bmod q)}}1
\EssBdd q\left( \sum_{0<|r|\le {q\over 2}}{1\over {|r|}}\right)\left({Q\over q}+1\right)\EssBdd Q
$$
\par
\noindent
(recall \thinspace $g\ge 0$, now), since \thinspace $(j,q)=1$ \enspace implies \enspace $jd\equiv r(\bmod \; q) \enspace \Leftrightarrow \enspace d\equiv r\overline{j}(\bmod \; q)$, where \enspace $\overline{j}(\bmod \; q)$ \enspace \thinspace is such that \enspace \thinspace  \thinspace $\overline{j}j\equiv 1(\bmod \; q)$, i.e., \enspace $\overline{j}$ \enspace is the {\stampatello reciprocal} of \enspace $j(\bmod \; q)$.
\par
See that we do not know the sign of \enspace $\Sigma_2(\alpha)$ \enspace (while we force \enspace $\Sigma_1(\alpha)\ge 0$ \enspace with \enspace $g\ge 0$), still we may estimate its \lq \lq maximum amplitude\rq \rq \thinspace and then make \enspace $S_f\ge 0$ \enspace {\stampatello simply choosing} \enspace $f(0)=QN^{2\varepsilon}$ \enspace (the number 
of \enspace $\varepsilon$s is unimportant, now):
$$
g\ge 0 \enspace \Rightarrow \enspace S_f(\alpha)= QN^{2\varepsilon} + \Sigma_1(\alpha) + \Sigma_2(\alpha) \ge 0 \quad \forall \alpha \in \hbox{\script F}_Q .
$$
\par
\noindent
Now, thanks to $(\ast)$, the result becomes (assuming above Theorem hypotheses, in particular \enspace $f\EssBdd 1$)
$$
g\ge 0 \enspace \Rightarrow \enspace I_f(N,h)\EssBdd I_{f,d_Q}(N,h)+Q^{2}h+Qh^{2}
$$
\par
\noindent
and resembles very much the results of [C1], though weaker (additional hypothesis \enspace $g\ge 0$ \enspace \& weaker error terms).
\par
\noindent
\centerline{\bf References}
\medskip

\item{\bf [C1]} \thinspace Coppola, G.\thinspace - \thinspace {\sl On the correlations, Selberg integral and symmetry of sieve functions in short intervals} \thinspace - \thinspace http://arxiv.org/abs/0709.3648v3 - 9pp. (electronic).
\smallskip

\item{\bf [C2]} \thinspace Coppola, G.\thinspace - \thinspace {\sl On the correlations, Selberg integral and symmetry of sieve functions in short intervals, II} \thinspace - \thinspace to appear.
\smallskip

\item{\bf [C-S]} \thinspace Coppola, G.\thinspace and \thinspace Salerno, S.\thinspace - \thinspace {\sl On the symmetry of the divisor function in almost all short intervals} \thinspace - \thinspace Acta Arith. {\bf 113} (2004), {\bf no.2}, 189--201. $\underline{\tt MR\enspace 2005a\!:\!11144}$\smallskip

\item{\bf [D]} \thinspace Davenport, H.\thinspace - \thinspace {\sl Multiplicative Number Theory} \thinspace - \thinspace Third Edition, GTM 74, Springer, New York, 2000. $\underline{{\tt MR\enspace 2001f\!:\!11001}}$
\smallskip

\item{\bf [T]} \thinspace Tenenbaum, G. \thinspace - \thinspace {\sl Introduction to analytic and probabilistic number theory} - Cambridge studies in advanced mathematics : 46, Cambridge University Press, 1995. $\underline{{\tt MR \enspace 97e\!:\!11005b}}$
\smallskip

\item{\bf [V]} \thinspace Vinogradov, I.M.\thinspace - \thinspace {\sl The Method of Trigonometrical Sums in the Theory of Numbers} - Interscience Publishers LTD, London, 1954. $\underline{{\tt MR \enspace 15,941b}}$
\medskip

\leftline{\tt Dr.Giovanni Coppola}
\leftline{\tt DIIMA - Universit\`a degli Studi di Salerno}
\leftline{\tt 84084 Fisciano (SA) - ITALY}
\leftline{\tt e-mail : gcoppola@diima.unisa.it}

\bye